\documentclass[a4paper,oneside,10pt]{amsproc} 

\usepackage[ansinew]{inputenc} 
\usepackage[T1]{fontenc}
\usepackage[english]{babel} 

\usepackage{enumerate} 
\usepackage{booktabs}

\usepackage{amsmath} 
\usepackage{amsfonts} 
\usepackage{amssymb}
\usepackage{amsthm}
\usepackage{latexsym}
\usepackage{bbm}
\usepackage{setspace}
\usepackage[top=3cm , bottom=3cm , left=3cm , right=3cm]{geometry}
\usepackage[pdftex]{graphicx}
  \DeclareGraphicsExtensions{.pdf,.png,.jpg,.jpeg,.mps}
  \usepackage{pgf}
  \usepackage{tikz}
\usepackage{floatflt}
\usepackage{url}
\usepackage{mathrsfs}
\usepackage{esint}
\usepackage{comment}

\makeatletter

\@namedef{subjclassname@2010}{%

  \textup{2010} Mathematics Subject Classification}

\makeatother

\theoremstyle{plain}
\newtheorem{lemma}{Lemma}
\newtheorem{theorem}{Theorem} 

\newtheorem*{prop*}{Proposition}
\newtheorem*{lemma*}{Lemma}

\newtheorem*{theorem*}{Theorem}

\theoremstyle{definition} 
\newtheorem*{def*}{Definition}

\theoremstyle{remark}
\newtheorem*{remark}{Remark}
\newtheorem*{acks}{Acknowledgements}

\newcommand{\R}{\mathbb{R}}

\newcommand{\Z}{\mathbb{Z}}

\newcommand{\E}{\mathbb{E}}

\newcommand{\D}{\,\textup{d}}

\newcommand{\Dd}{\mathscr{D}}

\newcommand{\Cc}{\mathscr{C}}
\newcommand{\Rr}{\mathscr{R}}

\newcommand{\Mm}{\mathscr{M}}
\newcommand{\la}{\langle}
\newcommand{\ra}{\rangle}

\newcommand{\BMO}{\textup{BMO}}
\newcommand{\loc}{\textup{loc}}
\newcommand{\supp}{\textup{supp}\,}

\title{Some remarks on the dyadic Rademacher maximal function}
\author{Mikko Kemppainen}
\address{Department of Mathematics and Statistics, University of Helsinki,
Gustaf Hällströmin katu 2b, FI-00014 Helsinki, Finland}
\email{mikko.k.kemppainen@helsinki.fi}

\begin{document}

\begin{abstract}
  Properties of a maximal function for vector-valued martingales were studied by the author in an earlier paper. Restricting here to the dyadic setting, we prove the equivalence between (weighted) $L^p$ inequalities and weak type estimates, and discuss an extension to the case of locally finite Borel measures on $\R^n$. In addition, to compensate for the lack of an $L^\infty$ inequality, we derive a suitable $\BMO$ estimate. Different dyadic systems in different dimensions are also considered.
\end{abstract}


\maketitle


\section{Introduction}

The \emph{Rademacher maximal function} was originally introduced by
Hytönen, McIntosh and Portal \cite{HMP} in order to prove a 
`Carleson's embedding theorem' for functions with values in infinite-dimensional Banach spaces. 
It provided a vector-valued analogue for the standard dyadic maximal function by replacing the 
suprema of local averages with
their R-bounds. 
More precisely, for locally integrable vector-valued functions $f$ on $\R^n$ they set
\begin{equation*}
  \Mm f(x) = \sup \Big\{ \Big( \E \Big\| \sum_{Q\ni x} \varepsilon_Q \lambda_Q \la f \ra_Q \Big\|^2
  \Big)^{1/2} : \Big( \sum_Q |\lambda_Q|^2 \Big)^{1/2} \leq 1 \Big\} , \quad x\in\R^n ,
\end{equation*}
where $\E$ denotes the expectation for independent random variables $\varepsilon_Q$ attaining 
values $+1$ and $-1$, each with probability $1/2$, and the vector $\la f \ra_Q$ is the average
of $f$ over a dyadic cube $Q\subset\R^n$.
The \emph{RMF property} of a Banach space $X$ was then defined by requiring that
for functions $f$ with values in $X$ we have
\begin{equation*}
  \int_{\R^n} \Mm f(x)^p \D x \lesssim\footnotemark \int_{\R^n} \| f(x) \|^p \D x ,
\end{equation*}
where $1 < p < \infty$.
\footnotetext{By $\alpha \lesssim \beta$ we mean that there exists a constant $C$ such that
$\alpha \leq C\beta$.
Quantities $\alpha$ and $\beta$ are comparable, $\alpha \eqsim \beta$, if
$\alpha \lesssim \beta$ and $\beta \lesssim \alpha$.}
In \cite{HMP} this property was shown to be independent of $p\in (1,\infty )$ 
and also to be non-trivial in the sense that while many spaces have it, 
not all (e.g. $\ell^1$) do.

The author studied this maximal function in a more general setting of martingales \cite{RMF} and showed,
employing somewhat lengthy arguments along the lines of \cite{MAUREYUMD} and \cite{BURKHOLDERUMD},
that the RMF property is characterized by a certain weak type estimate. A significantly
simpler approach is available if one restricts considerations to the original setting 
of dyadic cubes. 
Doing so enables us to extend the characterization of the RMF property and answer also other
natural questions concerning the Rademacher maximal function. Nevertheless, the question remains whether  
the RMF property follows from the better known UMD property --- a requirement for
unconditional convergence of Haar decompositions of vector-valued functions 
(see the remark on page \pageref{UMDdef}). 

Recently, the Rademacher maximal function has found applications in vector-valued $Tb$ theorems,
where one is typically led to study paraproduct operators, whose boundedness relies on
Carleson's embedding theorems. This was the case in an earlier version of 
\cite{HYTONENNONHOMTB} concerning a (global) vector-valued non-homogeneous $Tb$ theorem and in
a current version of its local counterpart \cite{HYTONENVAHLOCTB}.

The extended characterization of the RMF property is stated in Theorem \ref{RMF}, whereas 
Theorem \ref{BMO} entails the $\BMO$ estimate. Theorem \ref{systems} states the equivalence 
of $L^p$ inequalities with respect to different dyadic systems and the corresponding result for
different dimensions is presented in Theorem \ref{dimensions}. The characterization provided by 
Theorem \ref{RMF} is discussed in a more general setting of locally finite Borel measures in Section
\ref{gen}. 

\begin{acks}
  The financial support from Vilho, Yrjö and Kalle Väisälä Foundation is gratefully acknowledged.
  The author wishes to thank Tuomas Hytönen for helpful remarks.
\end{acks}


\subsection*{R-bounds}

Let $X$ be a Banach space and write $(\varepsilon_k)$ for a sequence of independent random variables 
attaining values $+1$ and $-1$, each with probability $1/2$. Comparison of randomized sums (and their expectations
$\E$) with square sums lies at the heart of our interest. 

\begin{def*}
A set $S\subset X$ is said to be \emph{R-bounded}\footnote{This coincides with the R-property in 
\cite[(2.4) Definition]{BERKSONGILLESPIE} when vectors $\xi$ are viewed
as operators $\lambda \mapsto \lambda\xi$ from scalars to $X$. The concept of R-boundedness appeared 
implicitly already in
\cite{BOURGAINDUALITY}.} if there exists a constant $C$ such that
\begin{equation*}
  \Big( \E \Big\| \sum_k \varepsilon_k \lambda_k \xi_k \Big\|^2 \Big)^{1/2} 
  \leq C \Big( \sum_k |\lambda_k |^2 \Big)^{1/2}
\end{equation*}
for all (finite) collections of vectors $\{ \xi_k \} \subset S$ 
and scalars $\{ \lambda_k \}$. The smallest such
$C$ is the R-bound $\Rr (S)$.
\end{def*}

\begin{remark}
\mbox{}
\begin{itemize}
\item
R-bounds satisfy the following `triangle inequality':
For $S,S'\subset X$ one has
\begin{equation*}
  | \Rr (S) - \Rr (S') | \leq \Rr (S\pm S') \leq \Rr (S) + \Rr (S') .
\end{equation*}
Furthermore, R-bounds are monotone and subadditive in the sense that
\begin{equation*}
  \Rr (S) \leq \sum_m \Rr (S_m) \quad \text{whenever} \quad S\subset \bigcup_m S_m \subset X .
\end{equation*}
In particular, for any sequence $( \xi_k )_{k=1}^\infty \subset X$ 
one has $\Rr ( \xi_1, \xi_2, \ldots ) \leq \sum_{k=1}^\infty \| \xi_k \|$.

\item 
R-bounds always exceed uniform bounds, that is,
\begin{equation*}
  \sup_{\xi\in S} \| \xi \| \leq \Rr (S) .
\end{equation*}
Moreover, that $\Rr (S) \lesssim \sup_{\xi\in S} \| \xi \|$ holds for all $S\subset X$ is equivalent with
$X$ having \emph{type} $2$ (see \cite[Proposition 1.13]{AB}). 
Recall that $X$ is said to have type $p\in [1,2]$ if
\begin{equation*}
  \Big( \E \Big\| \sum_k \varepsilon_k \xi_k \Big\|^2 \Big)^{1/2} 
  \lesssim \Big( \sum_k \| \xi_k \|^p \Big)^{1/p}
\end{equation*}
for all (finite) collections $\{ \xi_k \} \subset X$.

\end{itemize}
\end{remark}

\subsection*{The Rademacher maximal function}

Let us consider a system $\Dd = \bigcup_{k\in\Z} \Dd_k$ of \emph{dyadic cubes}, 
where each $\Dd_k$ partitions $\R^n$ into
cubes of sidelength $2^{-k}$ and every $Q\in\Dd_k$ is a union of $2^n$ smaller cubes $R\in\Dd_{k+1}$.
A standard example of such a system is given by $\Dd_k = \{ 2^{-k}([0,1)^n + m) : m\in\Z^n \}$. 
Note that every
$Q\in\Dd$ is contained in a unique larger cube $Q^*$ with $|Q^*| = 2^n|Q|$ ($| \cdot |$ 
refers to the Lebesgue measure)
and that for any two $Q,R\in\Dd$ the intersection
$Q\cap R$ is either $\emptyset$, $Q$ or $R$. By \emph{maximality} 
of a dyadic cube $Q$ in a given subcollection of $\Dd$
we mean that there does not exist a cube $R$ in the same subcollection for which $Q \subsetneq R$.
Note that maximal cubes are always disjoint and cover the same area as the whole subcollection.
Unless otherwise stated, $Q$ and $R$ will always stand for dyadic cubes in a given system.

For $1\leq p \leq \infty$, we denote by $L^p(X)$ the Lebesgue--Bochner space of $p$-integrable functions
(essentially bounded for $p=\infty$)
on $\R^n$ taking values in a Banach space $X$.

\begin{def*}
  The \emph{Rademacher maximal function} of an $f\in L^1_\loc (X)$ is given by
  \begin{equation*}
    \Mm f (x) = \Rr (\la f \ra_Q : Q\ni x ) , \quad x\in\R^n , \quad \text{where} \quad
    \la f \ra_Q = \frac{1}{|Q|} \int_Q f(y) \D y .
  \end{equation*}  
\end{def*}

\begin{remark}
\mbox{}
\begin{itemize}
\item
  If $X$ has type $2$, then R-bounds are comparable with uniform bounds and so $\Mm f$ is
  controlled pointwise by the standard dyadic maximal function
  \begin{equation*}
    Mf(x) = \sup_{Q\ni x} \| \la f \ra_Q \| .
  \end{equation*}
  
  \item \label{UMDdef} Functions with finite Haar decomposition form a dense subspace of $L^p(X)$ when 
  $1 < p < \infty$ and for such $f$ we have $\| \Mm f \|_{L^p} < \infty$. Recall, that every
  $f\in L^p(X)$ can be decomposed as
  \begin{equation*}
    f = \lim_{N\to\infty} 
    \sum_{\substack{Q\in\Dd_k \\ |k|\leq N}} \sum_\theta \la f,h_Q^\theta \ra h_Q^\theta ,
  \end{equation*}
  where the sum converges in $L^p(X)$ and the Haar functions $h_Q^\theta$ with 
  $\theta\in\{ 0,1 \}^n \setminus \{ 0 \}$ are defined as in
  \cite[Section 3]{VECTORVALUEDREVISITED}. In particular, each $h_Q^\theta$ is supported in $Q$, has
  $\int h_Q^\theta = 0$ and satisfies $|h_Q^\theta (x)| = |Q|^{-1/2}$ for all  $x \in Q$.
  Furthermore, $X$ is said to have the \emph{UMD property} if
  the convergence in the decomposition is unconditional in $L^p(X)$.
  
  \item Averages over large cubes have finite R-bounds for any $f\in L^p(X)$ with $1 \leq p < \infty$,
  that is, given any dyadic cube $Q$ we have $\Rr (\la f \ra_R : R\supset Q ) < \infty$.
  
  \item $\Mm$ preserves the dyadic support of functions with zero mean: If $\supp b \subset Q$ and
  $\int b = 0$, then for every $x\not\in Q$ and every $R\ni x$ we have $\la b \ra_R = 0$, since either
  $R\cap Q = \emptyset$ of $R\supset Q$.
  Consequently, $\Mm b(x) = 0$ for $x\not\in Q$.
  
\end{itemize}
\end{remark}

\section{$L^p$ inequalities and weak type estimates}

In this section we prove that for any Banach space $X$ and any 
$1 < p < \infty$, 
the $L^p$ inequality
\begin{equation*}
  \int_{\R^n} \Mm f(x)^p \D x \lesssim \int_{\R^n} \| f (x) \|^p \D x ,
\end{equation*} 
abbreviated as $\Mm : L^p(X) \to L^p$, is equivalent with
weak type estimates both on $L^1(X)$ and on the Hardy space $H^1(X)$. Moreover, we consider
weighted $L^p$ inequalities for weights in the (dyadic) Muckenhoupt classes $A_p$.

\subsection*{Weak type estimates}

The Hardy space $H^1(X)$ is taken to consist of those
$f\in L^1(X)$ for which the dyadic maximal function $Mf$ is integrable, so that the norm 
$\| f \|_{H^1(X)} := \| Mf \|_{L^1}$ is finite. An equivalent description is given in terms of 
\emph{atoms}:
A function $a\in L^q(X)$, where $1 < q \leq \infty$, is said to be a $q$-atom if there is a dyadic cube $Q$ so that
\begin{equation*}
  \supp a \subset Q , \quad \int_Q a(x) \D x = 0 , \quad \text{and} \quad \| a \|_{L^q(X)} 
  \leq |Q|^{-1/q'} ,
\end{equation*}
$q'$ being the Hölder conjugate of $q$. Note that every $q$-atom $a$ satisfies $\| a \|_{H^1(X)} \lesssim 1$. 
Now $H^1(X)$ consists of exactly those $f\in L^1(X)$ which admit, for every $q \in (1,\infty ]$, 
a decomposition into $q$-atoms $a_k$ so that
\begin{equation*}
  f = \sum_k \lambda_k a_k , \quad \text{with} \quad \sum_k |\lambda_k| < \infty .
\end{equation*}
The weak type Hardy space estimate is the requirement that
\begin{equation*}
| \{ x\in\R^n : \Mm f(x) > \lambda \} | \lesssim \frac{1}{\lambda} \| f \|_{H^1(X)} 
\end{equation*}
for all $\lambda > 0$. We write this as $\Mm : H^1(X) \to L^{1,\infty}$ (and similarly for $L^1(X)$).

The key to the derivation of an $L^p$ inequality from a weak type estimate 
is a suitable distributional inequality, where
$\Mm$ is controlled by another maximal operator.
For $1\leq q < \infty$ we define
\begin{equation*}
  M_q f(x) = \sup_{Q\ni x} \Big( \frac{1}{|Q|} \int_Q \| f(y) \|^q \D y \Big)^{1/q} .
\end{equation*}

\begin{lemma}
\label{distributional}
  Suppose that $\Mm : H^1(X) \to L^{1,\infty}$ and let $1 < q < \infty$.
  If $f$ has a finite Haar decomposition and $Q$ is maximal among cubes for which
  $\Rr (\la f \ra_R : R\supset Q ) > \lambda$ for a given $\lambda > 0$, then
  \begin{equation*}
    | \{ x\in Q : \Mm f(x) > 2\lambda , \: M_qf(x) \leq \delta \lambda \} |
      \lesssim \frac{\delta}{1-\delta} |Q|
  \end{equation*}
  for all $\delta \in (0,1)$. 
  Consequently, for every $\lambda > 0$ and $\delta\in (0,1)$, we have
  \begin{equation*}
    | \{ x\in\R^n : \Mm f(x) > 2\lambda , \: M_qf(x) \leq \delta \lambda \} | 
    \lesssim \frac{\delta}{1-\delta} | \{ x\in\R^n : \Mm f(x) > \lambda \} | .
  \end{equation*}
  \begin{proof}
    Given an $f$ with a finite Haar decomposition and a $\lambda > 0$,
    let $Q$ be maximal among cubes for which $\Rr (\la f \ra_R : R\supset Q ) > \lambda$. 
    
    If $\Mm f(x) > 2\lambda$ for an $x\in Q$, then $\Rr ( \la f \ra_R : R\subset Q , R\ni x ) > \lambda$,
    since $\Rr ( \la f \ra_R : R \supset Q^* ) \leq \lambda$ by maximality of $Q$. 
    If also $M_qf \leq \delta \lambda$ somewhere in $Q$, then
    \begin{align*}
      \Mm ((f-\la f \ra_Q )1_Q)(x) &= \Rr (\la f \ra_R - \la f \ra_Q : R\subset Q, R\ni x ) \\
      &\geq \Rr (\la f \ra_R : R\subset Q , R\ni x ) - \| \la f \ra_Q \| \\
      &> (1-\delta )\lambda ,
    \end{align*}
    as $\| \la f \ra_Q \| \leq M_qf(y)$ for any $y\in Q$.
   
    Now $(f - \la f \ra_Q)1_Q$ is $q$-atom multiplied by $2|Q|^{1/q'} \| f1_Q \|_{L^q(X)}$ and so
    from $\Mm :H^1(X) \to L^{1,\infty}$ it follows that
    \begin{align*}
      | \{ x\in Q : \Mm f(x) > 2\lambda , \: M_qf(x) \leq \delta \lambda \} |
      &\leq | \{ x\in Q : \Mm ((f-\la f \ra_Q )1_Q)(x) > (1-\delta )\lambda \} | \\
      &\lesssim \frac{1}{(1-\delta ) \lambda} \| (f-\la f \ra_Q )1_Q \|_{H^1(X)} \\
      &\lesssim \frac{1}{(1-\delta ) \lambda} |Q|^{1/q'} \| f1_Q \|_{L^q(X)} .
    \end{align*}
    Assuming that $M_qf \leq \delta \lambda$ somewhere in $Q$, we obtain
    \begin{equation*}
      \Big( \int_Q \| f(x) \|^q \D x \Big)^{1/q}
      \leq |Q|^{1/q} \inf_{x\in Q} M_qf(x) \leq |Q|^{1/q} \delta \lambda ,
    \end{equation*}
    so that from $|Q|^{1/q'}|Q|^{1/q} = |Q|$ we arrive at 
    \begin{equation*}
      | \{ x\in Q : \Mm f(x) > 2\lambda , \: M_qf(x) \leq \delta \lambda \} |
      \lesssim \frac{\delta}{1 - \delta} |Q| .
    \end{equation*}
    
    The set $\{ x\in\R^n : \Mm f(x) > \lambda \}$ can of course be decomposed into a disjoint union
    of maximal cubes in the previous sense and so
    \begin{equation*}
      | \{ x\in\R^n : \Mm f(x) > 2\lambda , \: M_qf(x) \leq \delta \lambda \} | 
      \lesssim \frac{\delta}{1-\delta} | \{ x\in\R^n : \Mm f(x) > \lambda \} |
    \end{equation*}
    for all $\delta \in (0,1)$.
  \end{proof}
\end{lemma}

\begin{remark}
  From $\Mm :L^1(X)\to L^{1,\infty}$ one can deduce a similar distributional inequality for $q=1$.
\end{remark}

\subsection*{Weights}

For $1 < p < \infty$, the (dyadic) Muckenhoupt class $A_p$ consists of \emph{weights} $w$ 
(non-negative and locally integrable) such that
\begin{equation*}
  \Big( \frac{1}{|Q|}\int_Q w(x) \D x \Big) \Big( \frac{1}{|Q|} \int_Q w(x)^{1-p'} \D x \Big)^{p-1} 
  \lesssim 1 
\end{equation*}
for every dyadic cube $Q$. This is equivalent to the requirement that, for any Banach space $X$,
$M_1: L^p(w;X) \to L^p(w)$, i.e.
\begin{equation*}
  \int_{\R^n} M_1f(x)^p w(x) \D x \lesssim \int_{\R^n} \| f(x) \|^p w(x) \D x .
\end{equation*}
Due to the `reverse Hölder property' of Muckenhoupt weights (see \cite[Chapter IV]{PEPEWEIGHTEDNORM}), every 
weight in $A_p$ belongs to a smaller class $A_{p/q}$ for some $q > 1$.
Furthermore, every such weight $w$ satisfies the following: 
There exists a $\gamma > 0$ such that, whenever $E\subset Q$ for a dyadic cube $Q$, we have
\begin{equation}
\tag{$\ast$}\label{fairshare}
  \frac{w(E)}{w(Q)} \lesssim \Big( \frac{|E|}{|Q|} \Big)^\gamma .
\end{equation}
Here, as usual, $w$ is also used to denote the measure $w(x)\D x$.

\subsection*{Characterization of the RMF property}
We are now in the position to characterize the RMF property of a Banach space by the equivalent
conditions in the following statement:

\begin{theorem}
\label{RMF}
  The following conditions are equivalent for any Banach space $X$:
  \begin{itemize}
  \item[(i)] $\Mm : L^p(w;X) \to L^p(w)$ for all $p\in (1,\infty )$ and any $w\in A_p$,
  \item[(ii)] $\Mm : L^p(X) \to L^p$ for some $p\in (1,\infty )$,
  \item[(iii)] $\Mm : L^1(X) \to L^{1,\infty}$,
  \item[(iv)] $\Mm : H^1(X) \to L^{1,\infty}$.
  \end{itemize}
  \begin{proof}
    As (ii) is a special case of (i), the equivalence is obtained by proving that 
    (ii) $\Rightarrow$ (iii) $\Rightarrow$ (iv) $\Rightarrow$ (i).
    
    (ii) $\Rightarrow$ (iii): To perform the Calderón--Zygmund decomposition for an
    $f\in L^1(X)$ at height $\lambda$, 
    let $\Cc$ denote the collection of maximal cubes among dyadic cubes $Q$ for which
    $1/|Q| \int_Q \| f (x) \| \D x > \lambda$. 
    We decompose $f$ into `good' and `bad' parts according to
    \begin{align*}
      g &= 1_{\R^n \setminus \bigcup\Cc} f + \sum_{Q\in\Cc} \la f \ra_Q 1_Q \\
      b &= f - g = \sum_{Q\in\Cc} (f- \la f \ra_Q)1_Q = \sum_{Q\in\Cc} b_Q .
    \end{align*}
     A standard argument employing the assumption $\Mm : L^p(X) \to L^p$ applies to the good part and gives
     \begin{equation*}
       | \{ x\in\R^n : \Mm g(x) > \lambda /2 \} | \lesssim \frac{1}{\lambda} \| f \|_{L^1(X)} .
     \end{equation*}
     For the bad part we observe that $\Mm b = 0$ outside $\bigcup\Cc$. Indeed,
     if $x\not\in\bigcup\Cc$ and $R\ni x$, then $\la b_Q \ra_R = 0$ for all $Q\in\Cc$ and so
     $\la b \ra_R = 0$.
     Consequently, also
     \begin{equation*}
       | \{ x\in\R^n : \Mm b(x) > \lambda /2 \} | \leq \Big| \bigcup\Cc \Big| 
       \leq \frac{1}{\lambda} \| f \|_{L^1(X)} .
     \end{equation*}
     
     (iii) $\Rightarrow$ (iv): This is immediate from the fact that 
     $\| \cdot \|_{L^1(X)} \leq \| \cdot \|_{H^1(X)}$.
     
     (iv) $\Rightarrow$ (i): Given a $p\in (1,\infty)$ and a $w\in A_p$, 
     we choose a $q \in (1,p)$ such that $w\in A_{p/q}$. 
     Any $f$ with a finite Haar decomposition will then satisfy,
     for all $\lambda > 0$ and $\delta \in (0,1)$, the inequality
     \begin{equation*}
       w( \{ x\in\R^n : \Mm f(x) > 2\lambda , \: M_qf(x) \leq \delta \lambda \} ) 
      \lesssim \Big( \frac{\delta}{1-\delta} \Big)^\gamma w( \{ x\in\R^n : \Mm f(x) > \lambda \} )  ,
     \end{equation*}
     with some $\gamma > 0$.
     Indeed, we may write $\{ x\in\R^n : \Mm f(x) > \lambda \}$ as a disjoint union of dyadic cubes $Q$
     that are maximal with respect to $\Rr (\la f \ra_R : R\supset Q ) > \lambda$, and
     then appeal to Lemma \ref{distributional} and to \eqref{fairshare} with
     $E = \{ x\in Q : \Mm f(x) > 2\lambda , \: M_qf(x) \leq \delta \lambda \}$ to see that
     there exists a $\gamma > 0$ so that
     \begin{equation*}
       w( \{ x\in Q : \Mm f(x) > 2\lambda , \: M_qf(x) \leq \delta \lambda \} )
       \lesssim \Big( \frac{\delta}{1-\delta} \Big)^\gamma w (Q) 
     \end{equation*}
     for all $\delta\in (0,1)$.     
     
     Now, writing $\alpha (\delta ) = (\delta / (1-\delta ))^\gamma$, we obtain
     \begin{align*}
       \| \Mm f \|_{L^p(w)}^p 
       &= 2^p \int_0^\infty p\lambda^{p-1} w( \{ x\in\R^n : \Mm f(x) > 2\lambda \} ) \D \lambda \\
       &\lesssim 2^p \alpha (\delta )
       \int_0^\infty p\lambda^{p-1} w( \{ x\in\R^n : \Mm f(x) > \lambda \} ) \D \lambda \\
       &\quad\quad 
       + 2^p \int_0^\infty p\lambda^{p-1} w( \{ x\in\R^n : M_qf(x) > \delta\lambda \} ) \D \lambda \\
       &= 2^p \alpha (\delta )
       \| \Mm f \|_{L^p(w)}^p + \frac{2^p}{\delta^p} \| M_qf \|_{L^p(w)}^p .
     \end{align*}
     Observing that $M_qf(x)^p = M_1g(x)^{p/q}$ for the scalar function $g(x)=\| f(x) \|^q$,    
     we may deduce from $w\in A_{p/q}$ that
     \begin{equation*}
       \| M_qf \|_{L^p(w)}^p = \int_{\R^n} M_1g(x)^{p/q} w(x) \D x
       \leq C_{p,q} \int_{\R^n} | g(x) |^{p/q} w(x) \D x = C_{p,q} \| f \|_{L^p(w;X)}^p.
     \end{equation*}      
     Choosing
     $\delta$ small enough so that $\alpha (\delta ) < 1/2^p$, we obtain
     after rearrangement that
     \begin{equation*}
       \| \Mm f \|_{L^p(w)}^p \lesssim \frac{2^p(C_{p,q})^p}{(1-2^p\alpha (\delta)) \delta^p} 
       \| f \|_{L^p(w;X)}^p .
     \end{equation*}
  \end{proof}
\end{theorem}

\begin{remark}
\mbox{}
\begin{itemize}
\item
  Condition (i) can also be seen to follow from (iii) 
  by using a distributional inequality as in Lemma \ref{distributional}, but with $q=1$.
  
  \item From condition (ii) it also follows that $\Mm : H^1(X) \to L^1$ as can easily be seen from the
  action of $\Mm$ on a $p$-atom $a$ supported in $Q$:
  \begin{equation*}
    \int_{\R^n} \Mm a(x) \D x \leq |Q|^{1/p'} \Big( \int_Q \Mm a(x)^p \D x \Big)^{1/p}
    \lesssim |Q|^{1/p'} \Big( \int_Q \| a(x) \|^p \D x \Big)^{1/p} \leq 1 .
  \end{equation*}
  
  \item The UMD property of a Banach space $X$ can be characterized by an analogous result for the 
  \emph{dyadic square function} given for $f\in L_\loc^1(X)$ by
  \begin{equation*}
    Sf(x) = \lim_{N\to\infty} \Big( \E \Big\| \sum_{\substack{Q\in\Dd_k \\ |k|\leq N}} \sum_\theta 
    \varepsilon_Q^\theta \la f,h_Q^\theta \ra h_Q^\theta (x)
    \Big\|^2 \Big)^{1/2}, \quad x\in\R^n ,
  \end{equation*}
  where $h_Q^\theta$ are the Haar functions:
  \begin{theorem*}
     The following conditions are equivalent for any Banach space $X$:
    \begin{itemize}
    \item[(i)] $S : L^p(w;X) \to L^p(w)$ for all $p\in (1,\infty )$ and any $w\in A_p$,
    \item[(ii)] $S : L^p(X) \to L^p$ for some $p\in (1,\infty )$,
    \item[(iii)] $S : L^1(X) \to L^{1,\infty}$,
    \item[(iv)] $S : H^1(X) \to L^{1,\infty}$.
    \end{itemize}
  \end{theorem*}
  The proof proceeds as that of Theorem \ref{RMF} once one has a suitable version of Lemma 
  \ref{distributional}. In order to prove a distributional inequality --- assuming that (iv) holds --- take 
  any $f$ with a finite Haar decomposition and a $\lambda > 0$. The set $\{ x\in\R^n : Sf(x) > \lambda \}$
  decomposes into disjoint cubes $Q$ that are maximal with respect to
  \begin{equation*}
    \Big( \E \Big\| \sum_{R\supset Q} \sum_\theta \varepsilon_R^\theta 
    \frac{\la f , h_R^\theta \ra }{|R|^{1/2}} \Big\|^2 \Big)^{1/2} > \lambda .
  \end{equation*}
  Now, if $Sf(x) > 2\lambda$ for an $x$ in such a cube $Q$, then
  \begin{align*}
    S((f-\la f \ra_Q)1_Q)(x) &= \Big( \E \Big\| \sum_{R\subset Q} \sum_\theta
    \varepsilon_R^\theta \la f , h_R^\theta \ra h_R^\theta (x) \Big\|^2 \Big)^{1/2} \\
    &\geq Sf(x) - \Big( \E \Big\| \sum_{R\supset Q^*} \sum_\theta \varepsilon_R^\theta 
    \frac{\la f , h_R^\theta \ra }{|R|^{1/2}} \Big\|^2 \Big)^{1/2} > \lambda ,
  \end{align*}
  where the identity follows from the fact that 
  $\la (f-\la f \ra_Q)1_Q , h_R^\theta \ra = \la f , h_R^\theta \ra$ if $R\subset Q$ and otherwise
  $\la (f-\la f \ra_Q)1_Q , h_R^\theta \ra = 0$.
  As in the proof of Lemma \ref{distributional}, we then have for any $q\in (1,\infty )$ that
  \begin{align*}
    |\{ x\in Q : Sf(x) > 2\lambda , \: M_qf(x) \leq \delta\lambda \}|
    &\leq |\{ x\in Q : S((f-\la f \ra_Q)1_Q)(x) > \lambda \}| \\
    &\lesssim \frac{1}{\lambda} \| (f-\la f \ra_Q)1_Q \|_{H^1(X)} \\
    &\lesssim \frac{1}{\lambda} |Q|^{1/q'} \| f1_Q \|_{L^q(X)} \\
    &\lesssim \delta |Q|,
  \end{align*}
  where the last step holds assuming that $M_qf \leq \delta\lambda$ somewhere in $Q$.
  
  Observe, in addition, that from $\| Sf \|_{L^p(w)} \lesssim \| f \|_{L^p(w;X)}$ one can deduce the
  reverse inequality $\| f \|_{L^p(w;X)} \lesssim \| Sf \|_{L^p(w)}$ by duality
  (cf. \cite[Theorem 5.4.7]{GRAFAKOS}).
\end{itemize}
\end{remark}

\subsection*{Application to paraproducts}
Let us briefly note how the weighted $L^p$ inequalities for $\Mm$ can be applied to 
vector-valued paraproducts.
We define the \emph{paraproduct operator} $\Pi_b$ associated to a given $b \in \BMO$ by
\begin{equation*}
  \Pi_b f = \sum_{Q,\theta} \la f \ra_Q \la b , h_Q^\theta \ra h_Q^\theta ,
\end{equation*}
where, strictly speaking, one considers finite sums and defines $\Pi_bf$ as a functional on a dense
subspace of the dual.
A standard argument via Carleson's embedding theorem (see \cite[Theorem 8.2, Corollary B.1]{HMP}
or \cite[Lemma 13, Theorem 14]{CARLESON}) gives
\begin{equation*}
  \| S(\Pi_bf) \|_{L^p(w)} \eqsim
  \Big( \E \Big\| \sum_{Q,\theta} \varepsilon_Q^\theta \la f \ra_Q \la b, h_Q^\theta \ra h_Q^\theta
  \Big\|_{L^p(w;X)}^p \Big)^{1/p} \lesssim \| b \|_{\BMO} \| \Mm f \|_{L^p(w)},
\end{equation*}
for $w\in A_p$ and $f\in L^p(w;X)$ with $1 < p < \infty$.
Assuming that $X$ has UMD, we have $\| \Pi_b f \|_{L^p(w;X)} \lesssim \| S(\Pi_bf) \|_{L^p(w)}$. 
If, in addition, $X$
has RMF, then $\| \Mm f \|_{L^p(w)} \lesssim \| f \|_{L^p(w;X)}$ according to Theorem \ref{RMF}, which
establishes the boundedness of $\Pi_b$ on $L^p(w;X)$. See \cite[Appendix B]{HMP} for historical remarks.

\section{A $\BMO$ estimate}

In contrast to other, more usual maximal operators (such as $M_q$), 
$\Mm$ does not in general map $L^\infty (X)$ boundedly into $L^\infty$. Indeed, according to 
\cite[Proposition 4.1]{RMF} we have:

\begin{prop*}
  For any Banach space $X$,
  $\Mm : L^\infty (X) \to L^\infty$ if and only if $X$ has type $2$.
\end{prop*}
On the other hand, a linearized version of $\Mm$ was shown in 
\cite[Proposition 7.1]{HMP} to map $L^\infty (X)$ into a 
certain vector-valued $\BMO$ space. Recall that by the John--Nirenberg inequality the (dyadic) $\BMO$ norm of an
$f\in L^1_\loc (X)$ can be given by any of the equivalent quantities
\begin{equation*}
  \| f \|_{\BMO (X)} \eqsim \sup_{Q\in\Dd} \Big( \frac{1}{|Q|} \int_Q \| f(x) - \la f \ra_Q \|^p \D x \Big)^{1/p} , \quad
  1 \leq p < \infty .
\end{equation*}
Moreover, the dyadic averages $\la g \ra_Q$ in the $\BMO$ norm of a scalar function $g\in L^1_\loc$
can be replaced by other scalars $c_Q$ according to the formula
\begin{equation*}
  \| g \|_\BMO \eqsim \sup_{Q\in\Dd} \inf_{c_Q} \frac{1}{|Q|} \int_Q | g(x) - c_Q | \D x .
\end{equation*}

\begin{theorem}
\label{BMO}
  Suppose that $\Mm : L^p(X) \to L^p$ for some $1 < p < \infty$. Then
  \begin{equation*}
    \| \Mm f \|_\BMO \lesssim \| f \|_{\BMO (X)}
  \end{equation*}
  for any $f \in L_\loc^1(X)$ with $\Mm f < \infty$ almost everywhere.
  \begin{proof}
    For every dyadic cube $Q$ and every $x\in Q$ we have
    \begin{align*}
      \Rr (\la f \ra_R : R\ni x) &\leq \Rr (\la f \ra_R - \la f \ra_Q + \la f \ra_{R'} : R\subset Q ,
      R\ni x , R'\supset Q) \\
      &\leq \Rr (\la f \ra_R - \la f \ra_Q : R\subset Q , R\ni x)
      + \Rr (\la f \ra_{R'} : R'\supset Q) ,
    \end{align*}
    where the first term in the last expression equals $\Mm ((f - \la f \ra_Q)1_Q)(x)$.
    Since $\Mm f < \infty$ almost everywhere, the constant
    \begin{equation*}
      c_Q = \Rr (\la f \ra_{R'} : R'\supset Q)
    \end{equation*}
    is finite and so for $x\in Q$,
    \begin{equation*}    
      0 \leq \Mm f(x) - c_Q \leq \Mm ((f - \la f \ra_Q)1_Q)(x) .
    \end{equation*}
    Consequently, since $\Mm : L^p(X) \to L^p$,
    \begin{align*}
      \frac{1}{|Q|} \int_Q |\Mm f(x) - c_Q| \D x 
      &\leq \frac{1}{|Q|} \int_Q \Mm ((f - \la f \ra_Q)1_Q)(x) \D x \\
      &\leq \Big( \frac{1}{|Q|} \int_Q \Mm ((f - \la f \ra_Q)1_Q)(x)^p \D x \Big)^{1/p} \\
      &\lesssim \Big( \frac{1}{|Q|} \int_Q \| f(x) - \la f \ra_Q \|^p \D x \Big)^{1/p}
      \lesssim \| f \|_{\BMO (X)} ,
    \end{align*}
    as required.
  \end{proof}
\end{theorem}


\section{Different dyadic systems and dimensions}
\label{sec:dim}

Up until now, we have considered the Rademacher maximal function with respect to a fixed
dyadic system on the Euclidean space of fixed dimension. 
It is shown in this section that the $L^p$ boundedness of $\Mm$ (as described in Theorem \ref{RMF}) 
depends neither on the system nor the dimension.

\subsection*{Different dyadic systems}
Different dyadic systems on $\R^n$ can be expressed by using a parameter
$\beta = (\beta_j) \in ( \{ 0,1 \}^n )^\Z$ according to $\Dd^\beta = \bigcup_{k\in\Z} \Dd^\beta_k$, with
\begin{equation*}
  \Dd^\beta_k = \{ 2^{-k}([0,1)^n +m) + \sum_{j>k} 2^{-j}\beta_j : m\in\Z^n \} .
\end{equation*}
The standard system corresponds to $\beta = 0$ and we refer to it by omitting $\beta$ altogether.
$L^p$ boundedness the Rademacher maximal operator with respect to $\Dd^\beta$ is
equivalent with uniform $L^p$ boundedness of the truncated operators defined by
\begin{equation*}
  \Mm^{\beta , N} f(x) = \Rr (A^\beta_kf(x) : k\geq N ) , \quad N\in\Z ,
\end{equation*}
where
\begin{equation*}
  A_k^\beta f = \sum_{Q\in\Dd^\beta_k} \la f \ra_Q1_Q
\end{equation*}
stands for an averaging operator with respect to $\Dd^\beta_k$. 
A direct calculation shows that averages with respect to $\Dd^\beta_k$ can be obtained from
those of the standard system by translations:
\begin{equation*}
  A^\beta_k = \tau^{-1}_k A_k \tau_k , \quad \text{where} \quad
  \tau_kf(x) = f(x + \sum_{j>k} 2^{-j}\beta_j ) .
\end{equation*}
Moreover, large (dyadic) translations commute with averaging so that for $k\geq j$ we have
\begin{equation*}
  A_k \sigma_j = \sigma_j A_k , \quad \text{where} \quad \sigma_jf(x) = f(x+2^{-j}\beta_j) .
\end{equation*}
Now that $\tau_{k-1} = \sigma_k \tau_k$, we see that actually
\begin{equation*}
  A^\beta_k = \tau^{-1}_N A_k \tau_N \quad \text{whenever} \quad k\geq N ,
\end{equation*}
and hence
\begin{equation*}
  \Mm^{\beta , N} f 
  = \tau^{-1}_N\Mm^N (\tau_Nf) .
\end{equation*}
Translations preserve $L^p$ norms and so we have arrived at the following result:


\begin{theorem}
\label{systems}
  Let $1 < p < \infty$.
  If the Rademacher maximal operator is $L^p$ bounded with respect to some dyadic system on $\R^n$,
  then it is $L^p$ bounded with respect to any dyadic system on $\R^n$.
\end{theorem}

\subsection*{Different dimensions}
Let us now consider the Rademacher maximal operator in different dimensions and prove the following result:

\begin{theorem}
\label{dimensions}
  Let $n$ be an positive integer and $1 < p < \infty$.
  The Rademacher maximal operator is $L^p$ bounded on $\R^n$ if and only if it is 
  $L^p$ bounded on $\R$ (or even on $[0,1)$).
\end{theorem}

We restrict our attention to the standard dyadic system on $\R^n$ and note that it
divides $\R^n$ into $2^n$ `quadrants' $\{ x\in\R^n : \alpha_j x_j \geq 0 \}$, where 
$\alpha\in \{ -1,1 \}^n$, in the sense that every cube in the standard system 
is contained in one of the (essentially disjoint) quadrants. For $L^p$ boundedness of $\Mm$ on $\R^n$,
it thus suffices to consider one of these quadrants, say $\{ x\in\R^n : x_j \geq 0 \}$.
By density of functions with bounded support, we may, using a scaling argument, restrict to functions
supported in the unit cube $[0,1)^n$ and consider only averages over cubes contained in $[0,1)^n$.
Writing $\Cc^n = \bigcup_{k=0}^\infty \Cc^n_k$, where $\Cc^n_k$ consists of (standard) dyadic cubes 
$Q\subset [0,1)^n$ of sidelength $2^{-k}$,
we have reduced the question to $L^p$ boundedness of
\begin{equation*}
  \Mm f(x) = \Rr ( \la f \ra_Q : Q\in\Cc^n , Q\ni x ) , \quad x \in [0,1)^n .
\end{equation*} 

To see that $\Mm$ is $L^p$ bounded on $[0,1)^n$ if and only if it is
$L^p$ bounded on $[0,1)$ we first note that `only if' is immediate from the fact that
functions on $[0,1)$ can be naturally viewed as functions on $[0,1)^n$ depending only on the first
coordinate. For sufficiency, we provide a way to associate dyadic subcubes of $[0,1)^n$ with
dyadic subintervals of $[0,1)$ in a suitable manner:

\begin{lemma}
\label{measpres}
  There exists a measure preserving map $\varphi : \Cc^n \to \Cc^1$ which respects the partial order
  of inclusions in the sense that for all $Q\in\Cc^n$, we have
  $\varphi (R) \subset \varphi (Q)$ if and only if $R\subset Q$. Moreover, for every $k\geq 0$,
  the restriction $\varphi_k : \Cc^n_k \to \Cc^1_{nk}$ is bijective.
  \begin{proof}
    Agreeing first that $\varphi ([0,1)^n) = [0,1)$, we proceed inductively. Namely,
    if $Q = 2^{-k}([0,1)^n + m)\in\Cc^n_k$ is mapped to 
    $\varphi (Q) = 2^{-nk}([0,1) + l)\in\Cc^1_{nk}$, 
    then each subcube $R\in\Cc^n_{k+1}$ of $Q$ is of the form
    \begin{equation*}
      R = 2^{-k-1}([0,1)^n + 2m + (\delta_1 , \ldots , \delta_n)) , 
      \quad \text{with} \quad \delta_j\in \{ 0,1 \},
    \end{equation*}
    and we map it to the interval
    \begin{equation*}
      \varphi (R) = 2^{-n(k+1)} ([0,1) + 2^nl + \delta_1 2^{n-1} + \cdots + \delta_n 2^0 ) ,
    \end{equation*}
    which is a subinterval of $\varphi (Q)$. Note that each subinterval $I\in\Cc^1_{n(k+1)}$
    of $\varphi (Q)$ is an image
    of exactly one subcube $R\in\Cc^n_{k+1}$ of $Q$ so that each restriction $\varphi_k$ is bijective.
  \end{proof}
\end{lemma}

Again, by switching to a truncation of $\Mm$, it suffices to consider, for each $N\geq 1$, 
functions on $[0,1)^n$ that are constant on cubes of $\Cc_N^n$.  
Every such $f$, when viewed as a function on cubes of $\Cc_N^n$,
can be transferred, using Lemma \ref{measpres}, 
to the function $f\circ \varphi_N^{-1}$ on $[0,1)$ (which is constant on cubes of
$\Cc^1_{nN}$). 
Dyadic averages of $f\circ \varphi_N^{-1}$ include the dyadic
averages of $f$; for every $Q\in\Cc_k^n$ with $0\leq k \leq N$ we have
\begin{equation*}
  \la f \ra_Q = \la f \circ \varphi_N^{-1} \ra_{\varphi (Q)} .
\end{equation*}
A calculation shows that the $L^p$ norm of $\Mm f$ is at most the $L^p$ norm of $\Mm (f\circ\varphi_N^{-1})$:
\begin{align*}
  \| \Mm (f\circ\varphi_N^{-1}) \|_{L^p([0,1))}^p 
  &= \frac{1}{2^{nN}} \sum_{J\in\Cc^1_{nN}} \Rr ( \la f \circ \varphi_N^{-1} \ra_I : I \supset J )^p \\
  &\geq \frac{1}{2^{nN}} \sum_{J\in\Cc^1_{nN}} 
     \Rr ( \la f \circ \varphi_N^{-1} \ra_{\varphi (Q)} : \varphi (Q) \supset J )^p \\
  &= \frac{1}{2^{nN}} \sum_{R\in\Cc^n_N} 
    \Rr ( \la f \circ \varphi_N^{-1} \ra_{\varphi (Q)} : \varphi (Q) \supset \varphi (R) )^p \\
  &= \frac{1}{2^{nN}} \sum_{R\in\Cc^n_N} 
    \Rr ( \la f \ra_Q : Q \supset R )^p = \| \Mm f \|_{L^p([0,1)^n)}^p .
\end{align*}
Since the $L^p$ norms of $f$ and $f\circ\varphi_N^{-1}$ are equal, Theorem \ref{dimensions} follows.

\section{More general measures}
\label{gen}

It was shown in \cite[Theorem 5.1]{RMF} that the RMF property of a Banach space $X$, as described here
by the equivalent conditions in Theorem \ref{RMF}, guarantees the boundedness of the Rademacher
maximal operator with respect to any filtration on any $\sigma$-finite measure space. It is
nevertheless interesting to see that the proof of Theorem \ref{RMF} is also directly applicable 
to a more general (possibly non-homogeneous) setting, where
$\R^n$ is equipped with a locally finite Borel measure $\mu$. We adjust our averages accordingly
by writing
\begin{equation*}
  \la f \ra_Q = \frac{1}{\mu (Q)} \int_Q f(y) \D\mu (y) , \quad Q\in\Dd ,
\end{equation*}
which we agree to be zero if $\mu (Q) = 0$, and put
\begin{equation*}
  A_kf = \sum_{Q\in\Dd_k} \la f \ra_Q 1_Q , \quad k\in\Z .
\end{equation*}


In order to obtain a collection of Haar functions $h_Q^\theta$ adapted to $\mu$, recall
that every $f\in L^p(\mu ; X)$ with $1 \leq p < \infty$ can be approximated by averages so that
$A_kf \to f$ in $L^p(\mu ; X)$ as $k\to\infty$.
For each $N\in\Z$ we can then write the truncated (adapted) Haar decomposition
\begin{equation*}
  f = \sum_{k=N}^\infty (A_{k+1}f - A_kf) + A_Nf
  = \sum_{\substack{Q\in\Dd_k \\ k\geq N}} \sum_\theta \la f , h_Q^\theta \ra h_Q^\theta
  + \sum_{Q\in\Dd_N} \la f \ra_Q1_Q ,
\end{equation*}
which converges (possibly conditionally) in $L^p(\mu ; X)$, 
and functions with finite decomposition are dense in 
$L^p(\mu ; X)$ (see \cite[Section 4]{HYTONENNONHOMTB}). 
Again, each $h_Q^\theta$ is supported in $Q$ and satisfies $\int h_Q^\theta \D\mu = 0$.

A suitable version of Lemma \ref{distributional}, with
\begin{equation*}
  M_1 f(x) = \sup_{Q\ni x} \frac{1}{\mu (Q)} \int_Q \| f(y) \| \D \mu (y) ,
\end{equation*}
can then be formulated as follows:

\begin{lemma}
\label{gendistributional}
  Suppose that $\Mm : L^1(\mu ; X) \to L^{1,\infty}(\mu)$ and let $N\in\Z$.
  If $f$ has a finite (adapted) Haar decomposition, then
  for every $\lambda > 0$ and $\delta\in (0,1)$, we have
  \begin{equation*}
    \mu ( \{ x\in\R^n : \Mm_N f(x) > 2\lambda , \: M_1f(x) \leq \delta \lambda \} ) 
    \lesssim \frac{\delta}{1-\delta} \, \mu ( \{ x\in\R^n : \Mm_N f(x) > \lambda \} ) ,
  \end{equation*}
  where $\Mm_Nf(x) = \Rr ( \la f \ra_Q : Q\ni x , Q\in\Dd_k, k\geq N )$ is the truncated Rademacher
  maximal function.
\end{lemma}

The truncation is needed in order to guarantee the existence of maximal cubes $Q$, which are now defined
by the requirement that $\Rr ( \la f \ra_R : R\supset Q, R\in\Dd_k, k\geq N ) > \lambda$.
Otherwise the proof proceed similarly to that of Lemma \ref{distributional}.


With these observations, the proof of Theorem \ref{RMF} can be adjusted to show the following generalization:

\begin{theorem}
\label{genRMF}
  Suppose that $\mu$ is a locally finite Borel measure on $\R^n$.
  The following conditions are equivalent for any Banach space $X$:
  \begin{itemize}
  \item[(i)] $\Mm : L^p(\mu ;X) \to L^p(\mu )$ for all $p\in (1,\infty )$,
  \item[(ii)] $\Mm : L^p(\mu ; X) \to L^p(\mu )$ for some $p\in (1,\infty )$,
  \item[(iii)] $\Mm : L^1(\mu ; X) \to L^{1,\infty}(\mu )$,
  \end{itemize}
  \begin{proof}
    (ii) $\Rightarrow$ (iii): We may argue as in \cite[Proposition 6.3]{RMF}.
    Given an $f\in L^1(\mu ; X)$ and a $\lambda > 0$, we show that
    \begin{equation*}
      \mu ( \{ x\in\R^n : \Mm_Nf(x) > \lambda \} ) \lesssim \frac{1}{\lambda} \| f \|_{L^1(\mu ; X)}
    \end{equation*}
    independently of $N\in\Z$. Gundy's decomposition (see \cite[Chapter IV, Section 2]{STEIN}
    or \cite[Theorem 6.2]{RMF}) allows us to write $f = g + h + b$, where
    \begin{enumerate}
      \item $\| g \|_{L^1(\mu ; X)} \lesssim \| f \|_{L^1(\mu ; X)}$ and 
      $\| g \|_{L^\infty (\mu ;X)} \lesssim \lambda$,
      \item $\| A_Nh \|_{L^1(\mu ; X)} + \sum_{k=N}^\infty
              \| A_{k+1}h - A_kh \|_{L^1(\mu ; X)} \lesssim \| f \|_{L^1(\mu ; X)}$,           
      \item $\mu ( \{ x\in\R^n : Mb(x) > 0 \} ) \lesssim \lambda^{-1} \| f \|_{L^1(\mu ; X)}$.
    \end{enumerate}
    From $\Mm : L^p(\mu ; X) \to L^p(\mu )$ and (1) it is straightforward to see that
    \begin{equation*}
      \mu ( \{ x\in\R^n : \Mm_Ng(x) > \lambda / 3 \} ) \lesssim \frac{1}{\lambda} \| f \|_{L^1(\mu ;X)}.
    \end{equation*}
    Also,    
    \begin{equation*}
      \mu ( \{ x\in\R^n : \Mm_Nb(x) > \lambda / 3 \} ) \lesssim \frac{1}{\lambda} \| f \|_{L^1(\mu ;X)}
    \end{equation*}
    follows immediately from (3) and the fact that $\Mm b (x) = 0$ if and only if $Mb(x) = 0$.
    
    In order to handle $\Mm_Nh$, we first observe that for any sequence of vectors 
    $(\xi_k)_{k=1}^\infty$ in $X$ we have
    \begin{equation*}
      \Rr \Big( \sum_{k=1}^j \xi_k : j\geq 1 \Big) \leq \sum_{k=1}^\infty \| \xi_k \| .
    \end{equation*}
    Thus for all $x\in\R^n$,
    \begin{equation*}
      \Mm_Nh(x) 
      \leq \| A_Nh(x) \| + \sum_{k=N}^\infty \| A_{k+1}h(x) - A_kh(x) \| ,
    \end{equation*}
    so that (2) gives
    \begin{equation*}
      \mu ( \{ x\in\R^n : \Mm_Nh(x) > \lambda / 3 \} ) \lesssim \frac{1}{\lambda} \Big(
      \| A_Nh \|_{L^1(\mu ; X)} + \sum_{k\geq N} \| A_{k+1}h - A_kh \|_{L^1(\mu ; X)}  \Big)
      \lesssim \frac{1}{\lambda} \| f \|_{L^1(\mu ; X)} .
    \end{equation*}
    Combining the estimates for $g$, $h$ and $b$ we obtain the desired result.
    
    (iii) $\Rightarrow$ (i): Given a $p\in (1,\infty )$ and an $N\in\Z$ we may use Lemma
    \ref{gendistributional} to see that for any $f$ with finite (adapted) Haar decomposition and
    any $\delta \in (0,1)$ we have
    \begin{equation*}
      \| \Mm_Nf \|_{L^p(\mu)}^p \lesssim 2^p \frac{\delta}{1 - \delta} \| \Mm_Nf \|_{L^p(\mu)}^p
      + \frac{2^p}{\delta^p} \| M_1f \|_{L^p(\mu)}^p ,
    \end{equation*}
    where both sides of the inequality are finite. Choosing $\delta$ small enough, we see that
    $\| \Mm_Nf \|_{L^p(\mu)} \lesssim \| f \|_{L^p(\mu ; X)}$ independently of $N$.
  \end{proof}
\end{theorem}

\bibliographystyle{plain}
\bibliography{viitteet}

\end{document}